\documentclass [12pt]{article}
\usepackage{latexsym,amsthm,amsfonts,amsmath}

\newtheorem{theo}{Theorem}[section]

\newtheorem{quest}[theo]{Question}
\newtheorem{lemma}[theo]{Lemma}

\newcommand\N{\ensuremath{\mathbb{N}}}
\newcommand\Z{\ensuremath{\mathbb{Z}}}

\newcommand{\ra}{\rightarrow}

\title{Almost Sure Recurrence of the Simple Random Walk Path}
 \author{Itai Benjamini \and Ori Gurel-Gurevich}     

\date{}

\begin{document}
\maketitle

\begin{abstract}
It is shown that the path of a simple random walk on any graph,
consisting of all vertices visited and edges crossed by the walk, is
almost surely a recurrent subgraph.
\end{abstract}


\section{Introduction}
Given a graph $G = (V,E)$ with finite degrees, a simple random
walk (SRW) on $G$ is a Markov chain on the set of vertices with
transition probabilities
$$
Prob( w_t = u | w_{t-1} = v) = 1/d_v,
$$
provided $\{u,v\} \in E$, where $d_v$ is the number of edges meeting
at $v$.

\bigskip
$G$ is called {\bf recurrent} iff a.s. SRW  visits any fixed vertex
infinitely often. It is called {\bf transient} otherwise.
\bigskip

Let $G$ be a graph. Let $PATH$ be the random subgraph of $G$,
consists of all vertices visited and edges crossed by a simple
random walk on $G$, that is, the path of the random walk.

\begin{theo}\label{main}
$PATH$ is a.s. recurrent.
\end{theo}

\begin{itemize}

\item For a recurrent $G$, the theorem is trivial, since any
subgraph of a recurrent graph is recurrent (see \cite{ds}). Also, in
that case $PATH=G$.

\item The theorem is already known for the Euclidean lattices,
since a.s. the SRW paths on three dimensional Euclidean lattice has
infinitely many cutpoints, i.e. points where the past of the path is
disjoint from its future, see \cite{ jp,l}. And then recurrence
follows by the Nash-Williams criterion \cite{nw}. An example of a
transient, bounded degree graph, for which $PATH$ has only finitely
many cutpoints a.s. is constructed in \cite{JLP,BGGS}.

\item
Morris \cite{mo} proved that the components of the Wired Spanning
Forest are a.s. recurrent, a result of similar spirit to the theorem
but with a different proof. For another a.s. recurrence theorem (for
distributional limits of finite planar graphs) see \cite{bs}.

\item
Exercise: show, without using theorem \ref{main}, that if $G$ is
transient then a.s. the SRW do not visit all the vertices of $G$.

\item
The proof uses the electrical networks interpretation of
recurrence.  For the connection between SRW and electrical network
see \cite{ds}. For further reading on recurrence see \cite{wo} and
the on-line lecture notes \cite{pe}.

\item One can think of a Brownian analogue of the theorem. That is
a.s. parabolicity of the Wiener sausage, with reflected boundary
conditions. It is of interest to formulate similar conjectures and
theorems for other generators and other random walks and
processes. For background on recurrence in the Riemannian context
see e.g. \cite{gr}.

For example, consider the range of a branching random walk on a
graph $G$, denoted by $R(BRW)$. Then we conjecture that almost
surely $R(BRW)$ is recurrent for BRW with the same branching law.
And a similar conjecture should hold for tree indexed random
walks. See \cite{bp} for definitions and background.

\end{itemize}

\begin{quest}
Given a graph $G$, denote by $PATH(n)$ the path created by the first
$n$ steps of the SRW on $G$. and by $R(n)$ the maximal electric
resistance between pairs of vertices on $PATH(n)$ (when $PATH(n)$ is
viewed as an electrical network where each edge is a one ohm
resistor).

By the theorem, on any bounded degree graph $R(n)\ra \infty$ a.s.
(note that $R(n) \ra \infty$ do not imply the theorem, e.g. balls in
the binary tree). Is there a uniform lower bound, depending on the
maximal degree, for the rate at which it grows, that is: Is there a
function $f$,
$$
\lim_n f(n) = \infty
$$

So that for any infinite graph of bounded degree. a.s.
$$
\limsup_n \frac{R(n)} {f(n)} > 0  ?
$$

\end{quest}

In particular one can speculate that $f(n) = C \log^2 n$ might work,
where the $\log^2 n$ is a lower bound coming from considering $R(n)$
when $G$ is $\Z^2$, which might be critical.

A different proof of theorem \ref{main} is provided in \cite{BGGL}.
In \cite{BGGS}, ideas from this paper and from \cite{BGGL} are
combined to provide some bounds on the resistance of the path on
finite segments of the graph.

The proof of the theorem is in the coming three sections. In the
next section we consider line-graphs with unbounded degrees.

\section{Proof of Theorem \ref{main} for line-graphs}

First, we shall prove the theorem for a very special case. Quite
surprisingly, the general case will not be very different.
Focusing on this special case will help illustrate the main ideas
of the proof.

A graph $G$ is called a \emph{line-graph} if $V_G=\N$ and $E_G$
includes only edges connecting successive vertices. Let $e_i$ denote
the number of edges connecting $i$ and $i+1$. We place {\bf no}
restriction on $e_i$.

\begin{theo}\label{special}
If $G$ is a line graph then $PATH$ on $G$ is a.s. recurrent.
\end{theo}
\begin{proof}
As always, the only interesting case is if $G$ is transient, which
is equivalent to $\sum_{i=0}^\infty {e_i}^{-1}<\infty$. Let $v(n)$
be the probability that a simple random walk starting at $n$
visits $0$. Clearly $v$ is a strictly decreasing function,
$v(0)=1$ and $lim_{n\rightarrow \infty} v(n) =0$. More precisely:

$$v(n)=\frac{\sum_{i=n}^\infty {e_i}^{-1}}
           {\sum_{i=0}^\infty {e_i}^{-1}}$$

$v$ is harmonic everywhere except at $0$. It follows that if $w_t$
is a simple (weighted) random walk, then the process $v(w_t)$ is
"almost" a martingale, i.e. it is a martingale as long as $w_t$ does
not reach $0$.

Let $s_n$ be the number of times the random walk crossed an edge
connecting $n$ and $n+1$, in either direction. Let $s'_n$ be the
number of edges connecting $n$ and $n+1$ which belong to $PATH$,
i.e. those edges that the random walk has crossed. The resistance of
$PATH$ is therefore $\sum_{i=0}^\infty {s'_i}^{-1}$. Obviously,
$s_n\ge s'_n$ so $\sum_{i=0}^\infty {s_i}^{-1}<\sum_{i=0}^\infty
{s'_i}^{-1}$. We will show that $\sum_{i=0}^\infty
{s_i}^{-1}=\infty$ almost surely, and therefore $PATH$ is almost
surely recurrent.

\begin{lemma}\label{special 0 1 law}
$Prob(\sum_{i=0}^\infty {s_i}^{-1}=\infty)$ is either $0$ or $1$.
\end{lemma}
\begin{proof}
Let $\{X_i^j\}_{i,j=0}^\infty$ be independent random variables,
defined by $Prob(X_i^j=1)=e_i/(e_{i-1}+e_i)$ and
$Prob(X_i^j=-1)=e_{i-1}/(e_{i-1}+e_i)$. Use these variables to
construct a simple random walk on $G$ in the obvious manner:
$w_{t+1}=w_t+X_{w_t}^t$. Now, $s_k$ is dependent (in the
probabilistic sense) only on $X_i^j$ for $i\ge k$, since every time
the walk is in $\{0,1,..,k-1\}$ it will almost surely reach $k$ at
some time. Therefore, a change to the values of finitely many of the
$X_i^j$s will change only the finitely many $s_i$'s and so cannot
effect the infiniteness of $\sum_{i=0}^\infty {s_i}^{-1}$. By
Kolomogorov's zero-one law we get that $Prob(\sum_{i=0}^\infty
{s_i}^{-1}=\infty)$ is either $0$ or $1$.
\end{proof}

It remains to show that $PATH$ is not almost surely transient.
First we shall handle the easy case, where the walk is quickly
transient.

\begin{lemma}\label{special simple case}
If for infinitely many $n$, $v(n)/2>v(n+1)$ then almost surely
$\sum_{i=0}^\infty {s_i}^{-1}=\infty$.
\end{lemma}
\begin{proof}
Let $\{n_i\}_{i=0}^\infty$ be an infinite series such that
$v(n_i)/2>v(n_i+1)$. Consider $p_i=Prob(s_{n_i}=1)$, the probability
that the random walk crosses an edge from $n_i$ to $n_i+1$ only
once. Let $\tau_i=\min(t|w_t=n_i+1)$ be the first time the random
walk reaches $n_i+1$. Let $\sigma_i=\min(t|t>\tau_i \cap w_t=n_i)$
be the first time after $\tau_i$ the walk reaches $n_i$ or $\infty$
if it never happens. Since $v$ is harmonic on $\{n_i,n_i+1..\}$ we
get that $\{v(w_t)\}_{t=\tau_i}^{\sigma_i}$ is a bounded martingale.
Adopting the convention $v(\infty)=0$, we get

$$v(n_i+1)=E(v(\tau_i))=E(v(\sigma_i))=0\cdot
Prob(\sigma_i=\infty)+v(n_i)\cdot Prob(\sigma_i<\infty)$$

Since $v(n_i+1)/v(n)<1/2$, the probability of ever reaching $n_i$
after having reached $n_i+1$ is less than $1/2$. This means that
$Prob(s_{n_i}=1)$ is at least $1/2$. By Fatou's lemma, the
probability of $s_{n_i}=1$ occurring infinitely often is at least
$1/2$ and so must be $1$ according to the proof of the previous
lemma. In particular, $\sum_{i=0}^\infty {s_i}^{-1}=\infty$ almost
surely.
\end{proof}

Lemma \ref{special simple case} shows that if $G$ is quickly
transient (in a rather weak sense) then $PATH$ almost surely has
infinitely many cut-edges and so must be recurrent.

If the premise of lemma \ref{special simple case} is not satisfied
then there must exist a sequence of vertices,
$\{n_i\}_{i=0}^\infty$, such that $n_0=0$ and
$v(n_i)/2>v(n_{i+1})>v(n_i)/4$.

Denote by $PATH_i$ the part of $PATH$ between $n_i$ and $n_{i+1}$.
Let $r_i=\sum_{j=n_i}^{n_{i+1}-1} {s_j}^{-1}$ be the resistance of
$PATH_i$.

Let
$$ q_i=\sum_{n_i\le w_t,w_{t+1}\le n_{i+1}} (v(w_{t+1})-v(w_t))^2 $$
i.e. the sum of $v(w_{t+1})-v(w_t))^2$ where the sum is taken over
the part of the random walk between $n_i$ and $n_{i+1}$.

Let $\tau_i=\min(t|w_t=n_i)$ be the first time the random walk
reaches $n_i$. Let $\sigma_i=min(t|t>t_i\cap w_t=n_{i-1})$ be the
first time after $t_i$ the random walk reaches $n_{i-1}$ or
$\infty$ if it never happens.

Let
$$ q'_i=\sum_{\tau_i \le t < \sigma_i} (v(w_{t+1})-v(w_t))^2 $$
i.e. the sum of $v(w_{t+1})-v(w_t))^2$ where the sum is taken over
the part of the random walk between times $\tau_i$ and $\sigma_i$.

\begin{lemma}\label{special expectation}
$$E(q'_i)<16 v^2(n_i)  $$
\end{lemma}
\begin{proof}
For prefixed $i$, let $a_t$ be equal to $v(w_{t+1})-v(w_t)$ if
$t<\sigma_i$ or 0 if $t\ge\sigma_i$. By definition
$v(n_i)+\sum_{t=\tau_i}^\infty a_t = v(w_{\sigma_i})$. Consider
$Var(v(w_{\sigma_i}))$. On the one hand we have

$$Var(v(w_{\sigma_i}))\le E(v^2(w_{\sigma_i})) \le v^2(n_{i-1}) < 16 v^2(n_i)$$

On the other hand

$$Var(v(w_{\sigma_i}))=\sum_{t=\tau_i}^\infty Var(a_t) +
2 \sum_{t=\tau_i}^\infty\sum_{t'=t+1}^\infty Cov(a_t,a_{t'})$$

By harmonicity of $v$, $E(a_t|w_0,w_1,..,w_t)=0$. Therefore
$Cov(a_t,a_{t'})=0$ for all $t\neq t'$.
$Var(a_t)=E((v(w_{t+1})-v(w_t))^2)$. Put together, we get

$$E(\sum_{\tau_i \le t < \sigma_i} (v(w_{t+1})-v(w_t))^2 =
\sum_{t=\tau_i}^\infty Var(a_t) = Var(v(w_{\sigma_i})) < 16
v^2(n_i)$$

\end{proof}

Now we use the connection between $q$ and $q'$ to prove the
following lemma.

\begin{lemma}\label{special lemma q}
$$Prob(q_i<64 v^2(n_i))>\frac{1}{4}$$
\end{lemma}
\begin{proof}
Using harmonicity of $v$ we get that
$Prob(\sigma_i<\infty)=v(n_i)/v(n_{i-1})<1/2$. From lemma
\ref{special expectation} we know that $E(q'_i)<16 v^2(n_i)$. $q'_i$
is nonnegative, so by Markov's inequality $Prob(q'_i<64
v^2(n_i))>3/4$. This implies that
$$Prob(\sigma_i=\infty \cap q'_i<64 v^2(n_i))>1/4$$

But if $\sigma_i$ is $\infty$ then $q'_i=q_i$ so
$$Prob(q_i<64 v^2(n_i))>1/4$$

\end{proof}

And finally we prove the relation between $q_i$ and $R_i$, the
resistance of $PATH_i$.

\begin{lemma}\label{special lemma R}
If $q_i<C v^2(n_i)$ then $R_i>\frac{1}{4C}$
\end{lemma}
\begin{proof}
Recall that $s_j$ is the number of times the walk crossed an edge
between $j$ and $j+1$. By definition
$$q_i=\sum_{j=n_i}^{n_{i+1}-1} s_j (v(j)-v(j+1))^2$$
and
$$R_i=\sum_{j=n_i}^{n_{i+1}-1} s_j^{-1} \ .$$

Using the Lagrange multipliers method, we try to minimize the value
of $R_i$, under the constraint given by the value of $q_i$. We get
$$\frac{\partial}{\partial s_j} ( R_i + \lambda q_i ) = -s_j^{-2}
+\lambda (v(j)-v(j+1))^2=0$$ which means that the minimum is
achieved when
$$s_j=\lambda^{-\frac{1}{2}}(v(j)-v(j+1))^{-1} \ .$$

Substituting $s_j$ in the definition of $q_i$ we get
$$q_i=\lambda^{-\frac{1}{2}} \sum_{j=n_i}^{n_{i+1}-1}
(v(j)-v(j+1))=\lambda^{-\frac{1}{2}}(v(n_i)-v(n_{i+1}))$$ which
implies
$$\lambda=(\frac{v(n_i)-v(n_{i+1})}{q_i})^2 \ .$$

Turning back to $R_i$ we get
$$R_i=\sum_{j=n_i}^{n_{i+1}-1} s_j^{-1} \ge \lambda^{\frac{1}{2}}\sum_{j=n_i}^{n_{i+1}-1}
(v(j)-v(j+1))=\frac{(v(n_i)-v(n_{i+1}))^2}{q_i}$$

$$>\frac{(v(n_i)-v(n_{i+1}))^2}{C v^2(n_i)}>\frac{v^2(n_i)}{4C
v^2(n_i)}=\frac{1}{4C}$$

\end{proof}

Now our work is nearly done. Combining lemma \ref{special lemma q}
and \ref{special lemma R} we get that for all $i$

$$Prob(R_i>\frac{1}{256})>\frac{1}{4}$$

Using Fatou's lemma again, we get

$$Prob(R_i>\frac{1}{256} \textsl{  infinitely often})>\frac{1}{4}$$

From lemma \ref{special 0 1 law} we know that the probability of
$PATH$ being recurrent is either $0$ or $1$. We just showed that it
cannot be 0 and therefore it must be $1$.

\end{proof}

\section{Proof of Theorem \ref{main} for bounded degree graphs}

Although the proof of theorem \ref{special} seems tailored to the
case of line graphs, only minor modifications are needed to adapt it
to the more general case of any bounded degree graph.

\begin{proof} 
First, we need to define $v$. Pick a vertex $g_0\in G$. Let $v(g)$
be the probability that a simple random walk starting at $g$ visits
$g_0$. For the general case it is not possible to give a simple,
closed formula for $v$, but it is easy to see that the relevant
properties of $v$ still hold: $v$ is harmonic except at $g_0$ and
$\lim_{t\rightarrow\infty}v(w_t)=0$ almost surely when $w$ is a
simple random walk.

Now we shall examine the four lemmas of the special case and prove
the corresponding lemmas for the general case.

Lemma \ref{special 0 1 law} proves a 0-1 law on the resistance of
$PATH$. While the conclusion of the lemma remain true for the
general case (we shall prove the resistance to be a.s. infinite),
the methods used in the proof are no longer valid. Indeed, it is not
true that the resistance of some part of $PATH$, far away from $g_0$
is a.s. independent of the "decisions" of the random walk made near
$g_0$. Instead of lemma \ref{special 0 1 law} we have the following
easy lemma.

\begin{lemma}\label{general 0 1 law}
If
$$Prob(PATH \textsl{ is transient})>0$$
then for every $C<1$ there exist a finite sequence of adjacent edges
$\overline{w}_0,..,\overline{w}_{t_0}$ such that
$$Prob(PATH \textsl{ is transient} \mid
(w_0,..,w_{t_0})=(\overline{w}_0,..,\overline{w}_{t_0}))>C$$
\end{lemma}
\begin{proof}
This is standard in measure theory. It follows easily from the
regularity of the random walk measure.
\end{proof}

Notice that all the arguments of the special case, as well as the
arguments we will use in the general case, can be carried out when
the random walk is conditioned to begin with a fixed sequence.

Next, we have lemma \ref{special simple case} which handles the
simple case where the walk is quickly transient. Here we don't have
this special case since we required the graph to have bounded
degree.

\begin{lemma}\label{general simple case}
If the degrees of vertices of $G$ are bounded by $d$, then for $g$
and $h$ adjacent vertices we have
$$v(h) \le d v(g)$$
\end{lemma}

\begin{proof}
This follows immediately from harmonicity of $v$.
\end{proof}

Let $C_i=\{g\in G \mid d^{-2i-1}\le v(g)\le d^{-2i}\}$ be the set of
all vertices whose $v$ values lies between $d^{-2i-1}$ and
$d^{-2i}$. From lemma \ref{general simple case} we know that every
$C_i$ is a cutset in the sense that it separates  $C_{i-1}$ from
$C_{i+1}$. It is not necessarily a cutset in the usual sense, of a
set separating $g_0$ from infinity, nor do these sets need be
finite. Indeed, there can be an infinite number of vertices for
which $v$ takes value above $d^{-2i}$. However, since $v(w_t)$ tends
to $0$ almost surely, the sets $C_i$ are cutset, in the usual sense,
in $PATH$ almost surely.

Let $PATH_i$ be all the edges in $PATH$ between $C_i$ and $C_{i+1}$.
More precisely,
$$PATH_i=\{(g,h)\in PATH \mid d^{-2i-2}<v(g)<d^{-2i-1} \ \bigcup \ d^{-2i-2}<v(h)<d^{-2i-1}\}$$

As before, let
$$q_i=\sum_{(w_t,w_{t+1})\in PATH_i} (v(w_{t+1})-v(w_t))^2$$
i.e. the sum of $v(w_{t+1})-v(w_t))^2$ over the part of the random
walk between $C_i$ and $C_{i+1}$.

Let $\tau_i=\min(t|w_t\in C_i)$ be the first time the random walk
reaches $C_i$. Let $\sigma_i=min(t|t>\tau_i\cap w_t\in C_{i-1})$ be
the first time after $\tau_i$ the random walk reaches $C_{i-1}$ or
$\infty$ if it never happens.

Let
$$ q'_i=\sum_{\tau_i \le t < \sigma_i} (v(w_{t+1})-v(w_t))^2 $$
i.e. the sum of $v(w_{t+1})-v(w_t))^2$ where the sum is taken over
the part of the random walk between times $\tau_i$ and $\sigma_i$.

\begin{lemma}\label{general expectation}
$$E(q'_i)< d^4 d^{-4i}$$
\end{lemma}
\begin{proof}
The proof is identical to that of lemma \ref{special expectation}.
This time we get
$$Var(v(w_{\sigma_i})) \le (d^{-2i+2})^2 = d^4d^{-4i}$$
and
$$Var(v(w_{\sigma_i})) = Var(v(w_{\tau_i}))+\sum_{t=\tau_i}^\infty
Var(v(w_{t+1})-v(w_t))$$
since the covariances are, as before, all
$0$.
\end{proof}

\begin{lemma}\label{general lemma q}
$$Prob(q_i< 4 d^4 d^{-4i}) \ge \frac{1}{4}$$
\end{lemma}
\begin{proof}
The proof is (again) identical to the proof of \ref{special lemma
q}. Here we have

$$Prob(\sigma_i<\infty) \le \frac{\sup_{g \in C_i} v(g)}{\inf_{g \in C_{i-1}} v(g)}\le \frac{1}{d} \le
\frac{1}{2}$$

and

$$Prob(q'_i<4 d^4 d^{-4i}) \ge \frac{3}{4}$$
\end{proof}

Next, we define $R_i$ as the resistance of $PATH_i$ when $C_i$ and
$C_{i+1}$ are both contracted, each to a single vertex, denoted
$c_i$ and $c_{i+1}$. The contracted $PATH_i$ will be denoted
$PATH'_i$.

\begin{lemma}\label{general lemma R}
If $q_i<C d^{-4i}$ then $R_i>\frac{1}{4Cd^2}$
\end{lemma}
\begin{proof}
The proof is actually simpler than \ref{special lemma R}. Let
$v'(g)$, defined for $g \in PATH_i$ be equal to $d^{-2i-1}$ for $g
\in C_i$, to $d^{-2i-2}$ for $g \in C_{i+1}$ and otherwise equal
to $v(g)$. By standard abuse of notation we shall refer to $v'$ as
defined on $PATH'_i$ too.

Let

$$q''_i=\sum_{(w_t,w_{t+1})\in PATH_i} (v'(w_{t+1})-v'(w_t))^2$$

Obviously, $q''_i \le q_i$. Now we use Thomson's Principle (see
\cite{ds} , page 49) on $PATH'_i$ with the function $v'$. $q''_i$ is
the "energy dissipation" of $v'$ on $PATH'_i$. By Thomson's
Principle the real energy dissipation is lower. Recall that
$v'(c_i)=d^{-2i-1}$ and $v'(c_{i+1})=d^{-2i-2}$.

Put together, we have
$$\frac{(d^{-2i-1}-d^{-2i-2})^2}{R_i} \le q''_i \le q_i < C d^{-4i}$$

Which yields
$$R_i>\frac{(d^{-2i-1}-d^{-2i-2})^2}{C d^{-4i}} \ge
\frac{1}{4Cd^2}$$
\end{proof}

Combining lemma \ref{general lemma q} and \ref{general lemma R} we
get that for all $i$
$$Prob(R_i>\frac{1}{16 d^6})>\frac{1}{4}$$

Using Fatou's lemma again we get that
$$Prob(R_i>\frac{1}{16 d^6} \textsl{  infinitely often})>\frac{1}{4}$$

By Rayleigh's Monotonicity Law (see \cite{ds} , page 51) we know
that the resistance of $PATH$ is greater than that of the
concatenation of $PATH'_i$, which is $\sum_{i=1}^\infty R_i$.
Therefore, the probability of $PATH$ being recurrent is greater
than $\frac{1}{4}$.

As noted earlier, all the arguments we used can be carried out
when the random walk is conditioned to begin with a fixed
sequence. Using lemma \ref{general 0 1 law}, we conclude that the
probability of $PATH$ not being recurrent must be $0$.

\end{proof}

\noindent {\bf Remark:} a close inspection of the proof reveals
that the theorem is also true for a finite union of paths of
independent simple random walks. The only difference is that lemma
\ref{general lemma q} applies to each SRW separately, to yield a
probability of $\frac{1}{4^k}$ ($k$ being the number of SRWs) for
the resistance of the union to be at least $\frac{1}{16 k d^6}$.

\section{Proof of Theorem \ref{main} general graphs}

Lastly, we turn our attention to the general case. Basically, what
happens here is that we forget about splitting our graph into
consecutive layers $C_i$, and instead consider the resistance
between just a single layer and infinity.

For this we need the following lemma:

\begin{lemma} \label{unbounded_resistance}
If $G$ is a transient graph then for every $\epsilon>0$ there is a
finite set of vertices $K_\epsilon$ such that the resistance between
$K_\epsilon$ and infinity is less then $\epsilon$.
\end{lemma}

\begin{proof}
Since $G$ is transient there is a unit flow from some vertex to
infinity with finite energy dissipation (see \cite{ds}, page 110).
Hence, there is a finite set $K_\epsilon$ such that the energy
dissipation outside of $K_\epsilon$ is less then $\epsilon$.
Thomson's principle implies that this bounds the resistance between
$K_\epsilon$ and infinity.
\end{proof}

Therefore, all we need is to show that with positive probability the
resistance between large balls in $PATH$ and infinity does not tend
to 0. For this it is enough to show that for uniformly for all balls
in $PATH$, the probability that the resistance to infinity is
bounded from below, is bounded from below.

So, let us define $v$ as in the previous section and let $w_t$ be a
simple random walk. Let $K$ be a finite set of vertices in $G$ and
let $v_K=\min(v(g)|g\in K)$ be the minimum voltage in $K$. define
$\overline{K}=\{g\in G | v(g)\ge v_K\}$ so that $\overline{K}\supset
K$. Let $R$ be the resistance, in $PATH$, between $\overline{K}$ and
infinity. Note that $\overline{K}$ might be infinite, but its
intersection with $PATH$ is a.s. finite.

Let $\tau=\min(t|v(w_t) < v_K)$ be the first time the walk exits
$\overline{K}$. Let $v_0=v(w_\tau)$. Let
$\sigma=\min(t|v(w_t)>2v_0)$, the first time the walk reaches twice
the voltage at $w_\tau$ or infinity if it never happens. $v(w_t)$ is
a martingale, so $Prob(\sigma=\infty)\ge \frac12$.

Given a vertex $g$ let $p_g=Prob(w_\sigma=g)$, so $\sum_g p_g
=Prob(\sigma<\infty)\le \frac12$. Let
$$q=\sum_{\tau\le t < \sigma} (v(w_{t+1})-v(w_t))^2$$

If $\sigma=\infty$ then in the definition of $q$ we sum over all
edges in $PATH$ that are not in $\overline{K}$. The following lemma
will bound the expectation of $q$ in that case.

\begin{lemma}
$$E(q|\sigma= \infty)\le 6 v^2_0$$
\end{lemma}

\begin{proof}
Defining, as before, $v(w_\infty)=0$ we get that
$$\sum_g p_g v(g) = E(v(w_\sigma))=v_0$$
and
$$Var(v(w_\sigma))= \sum_g p_g v^2(g) - v^2_0$$

On the other hand, the same argument as in lemma \ref{special
expectation} and \ref{general expectation} shows that
$$Var(v(w_\sigma))=\sum_{\tau\le t < \sigma} E((v(w_{t+1})-v(w_t))^2) = E(q)$$

Consider $E(q 1_{\sigma<\infty})$, where $1_{\sigma<\infty}$ is the
indicator function of the even $\sigma<\infty$. Obviously,
$$E(q 1_{\sigma=\infty})\ge E((v(w_\sigma)-v(w_{\sigma-1}))^2) \ge
E((v(w_\sigma)-2 v_0)^2)$$ since $v(w_{\sigma-1})\le 2 v_0$.
Expanding the expectation we get
$$E(q 1_{\sigma=\infty}) \ge \sum_g p_g (v(g) - 2 v_0)^2 = \sum_g
p_g v^2(g) - 4 v_0 \sum_g p_g v(g) + 4 v^2_0\sum_g p_g$$

$$=Var(v(w_\sigma))+v^2_0(4 Prob(\sigma < \infty) - 3)$$

Since
$$Var(v(w_\sigma))= E(q)= E(q1_{\sigma<infty}) +
E(q1_{\sigma=\infty})$$ we get that
$$E(q1_{\sigma=\infty})\le v^2_0(3-4 Prob(\sigma < \infty)) \le 3 v^2_0$$

Therefore,
$$E(q|\sigma=\infty)=\frac{E(q1_{\sigma=\infty})}{Prob(\sigma=\infty)}
\le 6 v^2_0$$
\end{proof}

Using this lemma we get that $Prob(q\le 12 v^2_0 | \sigma=\infty)
\ge \frac12$. Since $Prob(\sigma=\infty)\ge \frac12$, we know that
$Prob(q\le 12 v^2_0\ \cap \ \sigma=\infty)\ge \frac14$. If this
event happens then we have $R\ge \frac1{12}$, similarly to the proof
of lemma \ref{general lemma R}.

Using Fatou's lemma yet again, we have that the resistance in $PATH$
between balls and infinity does not tend to 0, with probability
$\frac14$. Applying lemma \ref{general 0 1 law} as in the previous
section concludes the proof.

\medskip
\noindent {\bf Acknowledgements:} Thanks to Gidi Amir, Gady Kozma,
Ron Peled and Benjy Weiss for useful discussions.

\end{document}